\documentclass[12pt]{article}
\usepackage[centertags]{amsmath}
\usepackage{amsfonts}
\usepackage{amssymb}
\usepackage{amsthm}
\usepackage{amsmath,mathrsfs}
\usepackage{amsmath,amscd}
\title{\textbf{Algebraic Stochastic Calculus}}
\author{Renaud Gauthier \footnote{rg.mathematics@gmail.com} \\ \\}
\theoremstyle{definition}

\newtheorem{Ff_elmt}{Proposition}[subsubsection]
\newtheorem{levelwise}[Ff_elmt]{Definition}
\newtheorem{filtGrothTop}[Ff_elmt]{Definition}
\newtheorem{filtsite}[Ff_elmt]{Definition}
\newtheorem{first_top}[Ff_elmt]{Proposition}
\newtheorem{ShBrown}{Definition}[subsection]
\newtheorem{Deform}[ShBrown]{Definition}
\newtheorem{cdga}{Definition}[subsection]
\newtheorem{top_proof}{Proposition}[subsection]
\newtheorem{HONO_top}{Proposition}[subsection]
\newcommand{\beq}{\begin{equation}}
\newcommand{\eeq}{\end{equation}}
\newcommand{\rarr}{\rightarrow}

\newcommand{\cA}{\mathcal{A}}
\newcommand{\cF}{\mathcal{F}}
\newcommand{\cO}{\mathcal{O}}
\newcommand{\cC}{\mathcal{C}}
\newcommand{\cG}{\mathcal{G}}
\newcommand{\cH}{\mathcal{H}}

\newcommand{\cD}{\mathcal{D}}
\newcommand{\cE}{\mathcal{E}}
\newcommand{\cR}{\mathcal{R}}
\newcommand{\cS}{\mathcal{S}}
\newcommand{\cW}{\mathcal{W}}
\newcommand{\opO}{\text{Op}(\Omega)}
\newcommand{\mR}{\mathbb{R}}
\newcommand{\mRf}{\mR^f}
\newcommand{\opOp}{\text{Op}(\Omega^+)}
\newcommand{\mapRcat}{\text{Map}(\mR, \text{Cat})}
\newcommand{\mapRopO}{\text{Map}(\mR, \opO)}
\newcommand{\mapRfopO}{\text{Map}(\mRf, \opO)}
\newcommand{\mapRfopOp}{\text{Map}(\mRf, \opOp)}
\newcommand{\mapRsA}{\text{Map}(\mR, \sA)}
\newcommand{\mapRsAO}{\text{Map}(\mR, \sAO)}

\newcommand{\mapRfsAO}{\text{Map}(\mRf, \sAO)}
\newcommand{\opORsA}{\opO^{\mR} \text{-} \sigma \text{alg}}
\newcommand{\opORfsA}{\opO^{\mRf} \text{-} \sA}
\newcommand{\opOpsdec}{\opOp \text{-} \sigma \text{dec}}
\newcommand{\opOpRfsdec}{\opOp^{\mRf} \text{-} \sigma \text{dec}}

\newcommand{\cov}{\text{Cov}}

\newcommand{\sA}{\sigma \text{-alg}}
\newcommand{\sAO}{\sA (\Omega)}
\begin{document}
\maketitle
\begin{abstract}
We develop the foundations of Algebraic Stochastic Calculus, with an aim to replacing what is typically referred to as Stochastic Calculus by a purely categorical version thereof. We first give a sheaf theoretic reinterpretation of Probability Theory. We regard probability spaces $(\Omega, \cF, P)$ as Grothendieck sites $(\cF, \tau_P)$ on which Brownian motions are defined via sheaves in symmetric monoidal $\infty$-categories. Due to the complex nature of such a formalism we are naturally led to considering a purely categorical, time independent formalism in which stochastic differential equations are replaced by studying problems in deformation theory.
\end{abstract}

\newpage

\section{Introduction}
Brownian motions display a fractal-like behavior (scale invariance (\cite{L})) and if one wants to use Brownian motions in modeling the dynamics of financial variables, this seems to contradict the fact that as time increments become smaller, the amount of information on the market one can collect should get smaller, and such information should be made more precise, hence a first motivation for departing from a purely Stochastic interpretation of Financial Mathematics. If financial instruments display a fractal-like behavior, one would want to resolve those phenomena. In order to do so one has to first investigate the formal foundations underlying the standard usage of financial variables. Our approach aims at explicating what it means to work with a probability space $(\Omega, \cF, P)$, what are Brownian motions $W$ defined on them, and how can one eliminate the sources of uncertainty, all of which contribute to a certain volatility, and from there eliminate the need for Brownian motions (also known as Wiener processes). Thus in a first time we deconstruct probability spaces $(\Omega, \cF, P)$. We regard $\Omega$ as the collection of objects (events) of a larger $\infty$-category $\Omega^+$ of events with $\cF$ seen as decorating objects of $\Omega$. We argue that $\cF$ can be regarded as a category on which we can put a Grothendieck topology, either operadic in nature, or induced by the probability $P$ itself, giving rise to a site $(\cF, \tau)$. Then we construct a sheaf-like object in the category of symmetric monoidal $\infty$-categories \cite{Lu1} on this site Brownian motions will be sections of. Due to the complex nature of such a sheaf, we consider instead a time independent formalism by working with $\Omega^+$ only with a Grothendieck topology on it induced by monomorphisms of simplicial sets. In this manner one can then construct a sheaf of symmetric monoidal $\infty$-categories whose sections do not have any volatility, thereby providing a categorical answer to the problem of having to work with Brownian motions. Alongside this development we also consider the tropical realization \cite{M} of the usual problem of having a geometric Brownian motion and study what that entails from the perspective of our formalism.\\

In the first section we deconstruct probability spaces $(\Omega, \cF, P)$ by regarding them as a set of objects $\Omega$ of some $\infty$-category $\Omega^+$ of events, which we decorate using a filtered $\sigma$-algebra $\cF = \{\cF_t\}$, and we regard the latter as an element of $\mapRcat$. We can put on $\cF$ two distinct Grothendieck topologies, one induced by focusing on the operadic aspect of $\cF$, the other by taking $P$ into consideration. If we denote by $\tau$ the latter Grothendieck topology, we then construct on $(\cF, \tau)$ a sheaf-like object $\cW$ in the $\infty$-category of $\infty$-categories on wich we put a symmetric monoidal structure, and which we regard as a sheaf of Brownian motions. This we do in section 3. In section 4, we abandon using $\sigma$-algebras and probabilities in favor of a purely categorical picture where $\tau_P$ is replaced by a Grothendieck topology induced by monomorphisms on events on a category of events $\cE$, and time is taken into consideration by considering roof diagrams under products $\omega \times \lhd(\omega)$, $\lhd(\omega)$ the forward cone at an event $\omega$ being the collection of events $\omega'$ that can be reached from $\omega$ via a morphism $\omega \rarr \omega'$. What we preserve from the Brownian motion picture however is the fact that for a suitable Grothendieck topology $\tau$ on $\cE$, we consider over $(\cE, \tau)$ sheaves of symmetric monoidal $\infty$-categories for which deformations of objects replace stochastic log differentials in stochastic calculus (\cite{KS1}, \cite{KS2}, \cite{O}).\\

\section{Algebraic Geometry and Probability Theory}
The starting point in this work is a random variable $X$ on a probability space $(\Omega, \cF, P)$. The log differential of such a variable sometimes depends on the differential of a Brownian motion. That latter is not well-defined however. We will consider the definition of such a differential and develop a formalism in which using such a differential takes its full meaning. In a first time we will be developing a sheaf theoretic formalism underlying the usual classical discussion of geometric Brownian motions, and that will serve as a platform towards a purely categorical approach to solving stochastic differential equations. This means we have to see probability spaces $(\Omega, \cF, P)$ as sites and we also have to define the sheaf over those sites Brownian motions are sections of.

\subsection{Probability space as a site}

\subsubsection{Filtered $\sigma$-algebras}

The set of all events $\Omega$ we regard as the set of objects of an underlying $\infty$-category of events. If the word ``event" is well-understood in a pedestrian sense, in our formalism we will model events as simplicial sets. The filtration $\cF$ we regard as organizing events according to a fashion dictated by the filtration itself. The filtration is a $\mR$-filtration. This gives an obvious map:
\begin{align}
\cF: \mR &\rarr \sA(\Omega) \nonumber \\
t &\mapsto \cF_t
\end{align}
where $\sA(\Omega)$ denotes the set of $\sigma$-algebras on a set $\Omega$. Thus $\cF \in \mapRsAO $. This filtration is increasing: $\cF_s \subseteq \cF_t$ for $s \leq t$. The manner in which information is collected however corresponds to summing graded sets; given $\cF_s$ for $s<t$, the information obtained at $t$ relative to $\lim_{\substack{s \rarr t \\ s < t}} \cF_s$ is really $\lim_{\substack{s \rarr t \\ s<t}} \cF_t / \cF_s$. \\

A filtration $\cF = \{\cF_t\}$ on $\Omega$ can be represented by an operad (\cite{MSS}, \cite{KM}) on $\Omega$, meaning the collection of information up to time $t$ can be mapped out using an operad. This we can do as we regard events to follow from an appropriate collection of previous events in the sense that for a given event $\omega$ there is at least one non-empty collection $\{ \omega_1, \cdots, \omega_N \}$ of events along with an appropriate map $\omega_1 \otimes \cdots \otimes \omega_N \rarr \omega$. We can assume $N$ to be finite; all the marginal contributions to $\omega$ from a possibly infinite collection of events can be repackaged into a single event $\omega_N$. Since events are indexed by $\mathbb{N}$ by having the index $t$ in $\cF_t$ being real we can by considering successive $\sigma$-algebras $\cF_t$ arrive at the above picture where $N$ is finite. Let $\opO$ denote the set of operads that can be defined on $\Omega$. Thus we have a map:
\begin{align}
\mathfrak{F}: \mapRopO & \rarr \mapRsAO \nonumber \\
(t \mapsto \cO(t)) &\mapsto (t \mapsto \cF_t)
\end{align}
where the map $\cO(t) \mapsto \cF_t$ gives $\cF_t$ the structure of a $\sigma$-algebra over an operad $\cO(t)$ through the collection of maps $\cO(t)(j) \otimes \cF_t ^j \rarr \cF_t$. Thus the value $\mathfrak{F}(\cO)$ of $\mathfrak{F}$ at an element $\cO$ of $\mapRopO$ gives the filtration $\cF$ on $\Omega$ the structure of a $\sigma$-algebra over an operad $\cO$. We regard $\cF$ as a $\mapRopO$-$\sigma$-algebra as there may be many operads $\cO$ that yield an $\cO$-algebra in $\mapRsAO$ that corresponds to $\cF$. Denote by $\opORsA$ the set of such filtered $\sigma$-algebras $\cF$ on $\Omega$.\\

\subsubsection{Nature of the filtration}
Time comes into the picture by virtue of having a filtration $\cF$ that is $\mR$-filtered. This means in particular that for a random variable $X$ on a probability space $(\Omega, \cF, P)$ for which $dX = \alpha X dt + \sigma X d W$, $dt$ in that equation is in reference to the index $t \in \mR$ in $\cF = \{\cF_t \}$. Further for a Brownian motion $W$ dependent on $t$ and at least one event $\omega$, $dW$ can be expressed in terms of $d \omega$, the deformation of $\omega$. Suppose $W$ depends on a unique event $\omega$. One can write $dW = adt + bd\omega$, where $dt$ is the same differential as in $dX$. Note that due to the Brownian motion nature of $W$, $a$ can be complicated but the nature of $a$ does not concern us here. However in Quantum Physics we have the well-known uncertainty relation $\Delta E \cdot \Delta t \geq \hbar /2$ (\cite{MT}). If energy translates into producing events, then $\Delta \omega \cdot \Delta t$ is bounded below, say by a constant $c>0$. Here we make the assumption that from a formal deformation $d \omega$ one can get a quantifiable increment $\Delta \omega$ by using an appropriately defined norm on $\Omega$. If $(dW)^2 = dt$, $(dt)^2 = 0$ as well as $dW dt = dt dW = 0$, as is typical to write in stochastic calculus (\cite{S}), $dt$ in $dW$ and the previous differential $dX$ cannot be the differential of physical time, otherwise $dW - adt = bd\omega$, $dt = b^2 d \omega \cdot d \omega$ and we would have $\Delta t \cdot \Delta \omega = b^2 (\Delta \omega)^3$ could be as small as we want. \\

Thus the indexing set $\mathbb{R}$ for the filtration does not correspond to a physical time, so the gathering of information dictates what is $t$ in $\cF_t$ as a function of physical time hence dictates how time must be contracted to obtain the indexing set for $\cF$. Events that continuously depend on physical time and are measured as such correspond to deformations of one given event, and events that do not continuously depend on physical time are measured on a discrete basis, so in all cases the collection of events corresponds to giving a discretization of the real line $\mathbb{R}$ corresponding to physical time, and if mod out that space by $\mathbb{N}$ we get a quotient map $\mR \rarr \mR/ \mathbb{N}$ with toroidal fibers, the fibers corresponding to intervals in time between collection times that come into the definition of individual $\sigma$-algebras $\cF_t$, $t \in \mR$. Here we are considering that collection times are discrete. This can be pictured geometrically by considering a framing on $\mR$ and accordingly we denote this quotient map by $\mR^f$. Note that if it is convenient to regard collection times as equivalent to having a quotient map $\mR \rarr \mR / \mathbb{N}$, it is far better geometrically to instead regard the situation as having a $(0,1]$-bundle over $\mR$, the base corresponding to the indexing set for $\cF$. We will still denote this bundle by $\mR^f$. To resolve the ambiguity in the manner in which events are collected together between consecutive collection times, one can extend $\cF$ to be defined on fibers of $\mR^f$ hence we define the extension $\cF^f$ of $\cF$ from a  $\mapRopO$-$\sigma$-algebra to a $\mapRfopO$-$\sigma$-algebra. Denote by $\opORfsA$ the set of such algebras. Thus $q: \mR^f \rarr \mR$ induces:
\begin{align}
q^*: \mapRsAO & \rarr \mapRfsAO \nonumber \\
\cF & \mapsto q^*(\cF)
\end{align}
and more generally:
\beq
\cF^f: \mapRfopO \rarr \mapRsAO \label{eff}
\eeq
is the extension of $\cF$ whose morphisms between events in the fibers of $\mR^f \rarr \mR$ are dictated by filtered operads in $\mapRfopO$.\\

\subsubsection{Probabilities on filtered $\sigma$-algebras}

Note that working with $(\Omega, \cF, P)$ is really the same as defining $P_t$ on $\cF^f (-, t)=\cF^f _t$ and:
\beq
P = \lim_{\rarr}P_t
\eeq
$P$ as used in Finance is a sub-homomorphism:
\beq
P: (\cF_{\infty}^f, \pi \circ \text{mult}_{\cO}) \rarr ([0,1], +)
\eeq
where $\text{mult}_{\cO}$ is the multiplication on $\cO$, an operad that defines $\cF^f_{\infty}$ and $\pi$ is the forgetful map that projects the $\cO$-$\sigma$-algebra $\cF^f_{\infty}$ to its underlying set. We have for distinct $\omega_i \in \cF^f_{\infty}$, $1 \leq i \leq j$:
\beq
P(\pi \circ \text{mult}_{\cO}(\omega_1, \cdots , \omega_j)) = P(\coprod_i \omega_i) =  \sum_i P(\omega_i)
\eeq
For generic $\omega_i$'s however:
\beq
P(\pi \circ \text{mult}_{\cO}(\omega_1, \cdots , \omega_j)) = P(\cup_i \omega_i) \leq  \sum_i P(\omega_i)
\eeq
If we denote $\text{Sub-Hom}(\cF^f , [0,1])$ by $(\cF^f)^{\underline{\wedge}}$, then we can represent this formalism with a diagram:
\beq
\begin{CD}
\cF^f _t @>\text{colim} >> \cF^f_{\infty} \\
@V\wedge VV @VV\wedge V \\
P_t \in (\cF^f _t )^{\underline{\wedge}} @>>> (\cF^f_{\infty})^{\underline{\wedge}} \ni P
\end{CD}
\eeq
Then a real valued random variable $X$ is an element of $\text{Mor}((\cF_{\infty}^f \times (\cF_{\infty}^f)^{\underline{\wedge}}), \mR) = (\cF_{\infty}^f \times (\cF_{\infty}^f)^{\underline{\wedge}})^{\wedge}$.\\

Now $P$ is really defined by its value on connected components of the $\opO ^{\mR^f}$-$\sigma$-algebra $\cF^f _{\infty}$. Indeed, two events that are unrelated in a loose sense cannot be compared, hence for two such events $\omega$ and $\omega'$ saying $P(\omega) < P(\omega')$ is not helpful. However for two events $\omega$ and $\omega'$ that are elements of a same connected component saying $P(\omega) < P(\omega')$ does correspond to a useful statement. This means in particular $P$ is not defined on the set underlying $\cF^f_{\infty}$ but on the $\opO^{\mR^f}$-$\sigma$-algebra $\cF^f_{\infty}$ itself on which one has to define a notion of connected components. Keeping the operadic picture in mind we regard $\Omega$ as the set of objects of an $\infty$-category $\Omega^+$ of events and higher morphisms between them. Further if we denote by $\cE v$ the category of $\infty$-categories of events, we take $\Omega^+$ to be a universal object in that category in the sense that $\text{Ob}(\Omega^+)$ is the set of all events and for any two events $\omega$ and $\omega'$ if there is a map $f: \omega \rarr \omega'$ then $f \in \text{Mor}_{\Omega^+}(\omega, \omega')$. This leads to defining an extension of $\cF^f$ in the following sense. An element $\cF$ of $\mapRsAO$ is now seen as a decoration of objects of $\Omega^+$ via a map $\mapRsA \rarr \Omega^+$. Such maps form what we denote by $\sigma \text{dec}(\Omega)$ and which we refer to as the set of $\sigma$-decorations of $\Omega$. This carries over to the set of $\cF^f$'s. Each element of that set is seen as a decoration over an operad in $\mapRfopOp$ and the set of such elements we denote by $\mapRfopOp \text{-}\sigma \text{dec}$ or $\opOpRfsdec$ for short.\\

We have:
\beq
P \in (\opOpRfsdec)^{\underline{\wedge}}
\eeq
Finally a Brownian motion $W$ is an element of $(\opOpRfsdec \times_{\opOpRfsdec} \opOpRfsdec^{\underline{\wedge}})^{\wedge}$. If we tentatively let $\cH$ be the space in which Brownian motions are valued, then observe that $W_*P = P(W^{-1})$ defines a push-forward measure on $\cH$.\\

\subsubsection{Filtered $\sigma$-algebras as sites}

Now fix some $\cF^f \in \opOpRfsdec$.
\begin{Ff_elmt}
$\cF^f \in \opOpRfsdec$ is an element of $\mapRcat$.
\end{Ff_elmt}
\begin{proof}
For each $t \geq 0$, $\cF^f_t$ provides a decoration of Ob($\Omega^+$) thereby defining a category $\cF^f_t$ whose objects are $\cF^f_t$-decorated elements of $\text{Ob}(\Omega^+)$, morphisms are morphisms of $\Omega^+$ between $\cF^f_t$-decorated elements if any, $\text{id} = \text{id}_{\Omega^+}$, composition is associative, the identity is a right and left unit, and this for all $t \geq 0$, making $\cF^f = \{ \cF^f_t \}_{t \geq 0}$ an element of $\text{Map}(\mR_{ \geq 0}, \text{Cat})$.
\end{proof}
\begin{levelwise}
We define connected components of $\cF^f$ to be collections of levelwise $\cF^f$-decorated connected components of $\Omega^+$, $t \geq 0$, of the form $\{ \{C(t)\} \quad | \quad t \geq 0  \}$ where $C(t)$ is a $\cF^f _t$-decoration of a connected component $C$ of $\Omega^+$.
\end{levelwise}
\begin{filtGrothTop}
For $S$ a set, $f \in \text{Map}(S, \text{Cat})$, we call a collection $\tau = \{\tau(u) \}_{u \in S}$, $\tau(u)$ a Grothendieck topology on $f(u) \in \text{Cat}$ for all $u \in S$, a filtered Grothendieck topology $\tau$ on $f$.
\end{filtGrothTop}

We will now define a filtered Grothendieck topology on $\cF^f$ that ties in the underlying $\infty$-category $\Omega^+$ with the operadic aspect of $\cF^f$.\\

Levelwise products in a filtered operad defining $\cF^f$ define an obvious filtered Grothendieck topology $\tau$ on $\cF^f$ as we will see. We define t-coverings of an event $\omega$ to be given by a collection $\cov^{(t)}(\omega) = \{\omega' \rarr \omega \}$ where $\omega'$ is in the same connected component $C(t)$ of $\cF_t^f$ $\omega$ is in and there is a multiplication in an operad generating $\cF^f_t$ one of whose summands is $\omega'$ and whose target is $\omega$. Such coverings will define a Grothendieck topology $\tau(t)$ the collection of which assembles into a filtered Grothendieck topology $\tau = \{ \tau(t) \}_{t \geq 0}$. \\

\begin{filtsite}
For $S$ a set, $f \in \text{Map}(S, \text{Cat})$, $\tau$ a filtered Grothendieck topology on $f$, we call $(f, \tau)$ a filtered site.
\end{filtsite}

\begin{first_top}
$\tau$ defines a filtered Grothendieck topology on $\cF^f \in \opOpRfsdec$ making $(\cF^f, \tau)$ into a filtered site.
\end{first_top}
\begin{proof}
If $\omega' \rarr \omega$ is an isomorphism then $\{\omega' \rarr \omega \}$ is a t-covering for all $t$. If $\{\omega_i \rarr \omega \}$ is a t-covering, $\gamma \rarr \omega$ is any arrow in $\cF^f_t$, by the universality of $\Omega^+$ since $\omega_i \times_{\omega} \gamma$ is an event there is a morphism $\omega_i \times_{\omega} \gamma \rarr \gamma$ in $\Omega^+$. It suffices now to show $\omega_i \times_{\omega} \gamma \in \cF_t ^f$ and there is a product in an operad defining $\cF_t ^f$ sending $\omega_i \times_{\omega} \gamma$ to $\gamma$. Since $\cF_t^f$ is a $\sigma$-algebra if $\omega_i$ and $\gamma$ are in $\cF^f_t$ so is $\omega_i \cup \gamma$. If $\omega_i$ and $\gamma$ are distinct then $\omega_i \cup \gamma \simeq \omega_i \times \gamma \in \cF_t^f$, assembled at that point. In particular there is at some time $s \leq t$ a morphism $\omega_i \times_{\omega} \gamma \rarr \omega_i \times \gamma$, so $\omega_i \times_{\omega} \gamma \in \cF_s^f$, in particular it is in $\cF^f_t$. If $\omega_i = \gamma$ however, $\omega_i \times_{\omega} \omega_i = \omega_i \times \omega_i$, which as an object is equivalent to considering $\omega_i$. In all cases then $\omega_i \times_{\omega} \gamma \in \cF^f_t$. Since we have maps $\omega_i \rarr \omega$ and $\gamma \rarr \omega$, we have an operad giving a map $\omega_i \otimes \gamma \rarr \omega$, and a universal such map we can find in an operad $\cO(s)$ for $s \leq t$ with a commutative diagram:
\beq
\begin{CD}
\omega_i \times_{\omega} \gamma @>>> \gamma \\
@VVV @VVV \\
\omega_i @>>> \omega
\end{CD}
\eeq
The top map $\omega_i \times_{\omega} \gamma \rarr \gamma$ in in $\cO(s)$, so it yields a morphism in $\cF^f_s$, in particular it is in $\cF^f_t$, and consequently we have an operad $\cO(t)$ giving rise to a map $\omega_i \times_{\omega} \gamma \rarr \gamma$. Thus $\{ \omega_i \times_{\omega} \gamma \rarr \gamma \}$ is an element of $\cov^{(t)}(\gamma)$. Finally if $\{\omega_i \rarr \omega\}$ is a covering $\cov^{(t)}(\omega)$, $\{\omega_{ij} \rarr \omega_i \}$ is a covering $\cov^{(t)}(\omega_i)$ then $\{\omega_{ij} \rarr \omega_i \rarr \omega\}$ is a covering $\cov^{(t)}(\omega)$. Thus $\tau(t)$ defines a Grothendieck topology on $\cF^f_t$, and this for all $t$, hence $\tau = \{ \tau(t) \}_{ t \geq 0}$ is a filtered Grothendieck topology on $\cF^f = \{ \cF^f_t \}_{t \geq 0}$.
\end{proof}

\subsection{Probability space as a site - version 2}
What we have achieved in the previous subsection is really get an enhanced site, a category $\cF^f$ with a Grothendieck topology and a probability $P$ defined on it as well. In the present subsection we use $P$ to define a Grothendieck topology $\tau_P$ on $\cF^f$ instead of trying to tie in the operadic picture. This has the advantage of working with $(\cF^f, \tau_P)$ instead of $(\cF^f \times (\cF^f)^{\wedge})$ as it is implied in Finance. The Grothendieck topology $\tau_P$ on $\cF^f \in \opOpRfsdec$ induced by $P$ is defined as follows. t-coverings of $\omega \in \cF^f$ are defined levelwise by $\cov^{(t)}(\omega) = \{\omega' \rarr \omega \, | \, \omega' \in C(t), P(\omega) \geq P(\omega')\}$ with $\omega \in C(t)$. We will first show those define a Grothendieck topology on $\cF^f_t$, and then we will define $\tau_P = \{ \tau_P(t) \}_{ t \geq 0}$ as giving a filtered Grothendieck topology on $\cF^f = \{ \cF^f_t \}_{ t \geq 0} \in \mapRcat$. \\

\begin{top_proof}
$\tau_P$ defines a filtered Grothendieck topology on a given $\cF^f \in \opOpRfsdec$ making $(\cF^f, \tau_P)$ into a filtered site.
\end{top_proof}
\begin{proof}
If $\omega' \rarr \omega$ is an isomorphism, then $P(\omega') = P(\omega)$ so $\{\omega' \rarr \omega \}$ is a covering for all $t$. If $\{\omega_i \rarr \omega \}$ is a t-covering, $\gamma \rarr \omega$ any arrow in $\cF^f_t$, since $P(\omega_i \times_{\omega} \gamma) \leq P(\omega_i \times \gamma) \leq P(\gamma)$, $\omega_i \times_{\omega} \gamma \rarr \gamma$ is a morphism in $\Omega^+$, $\omega_i \times_{\omega} \gamma$ is an element of $\cF^f_t$, then $\omega_i \times_{\omega} \gamma \rarr \gamma$ is a morphism in $\cF^f_t$ as observed before, so that $\{ \omega_i \times_{\omega} \gamma \rarr \gamma \}$ is a t-covering $\cov^{(t)}(\gamma)$. Finally if $\{\omega_i \rarr \omega \}$ is a t-covering and so is $\{ \omega_{ij} \rarr \omega_i \}$ then the morphism $\omega_{ij} \rarr \omega_i \rarr \omega$ is a morphism in $\cF^f_t$ and $P(\omega_{ij}) \leq P(\omega_i) \leq P(\omega)$ so that $P(\omega_{ij}) \leq P(\omega)$ and thus $\{\omega_{ij} \rarr \omega_i \rarr \omega \}$ is a t-covering, hence we have a Grothendieck topology $\tau_P(t)$ on $\cF^f_t$, and this for all $t$, so we can assemble the $\tau_P(t)$ into $\tau_P = \{ \tau_P(t) \}_{ t \geq 0}$, a filtered Grothendieck topology on $\cF^f$.
\end{proof}

We will denote $(\cF^f, \tau_P)$ by $(\cC, \tau)$ for simplicity and the probability will be implied. This will be preferred over working with $(\cF^f, \tau)$ from the previous subsection which is actually the object being implied in Quantitative Finance.\\

\subsection{Time independent formalism} \label{TIF}
As we will argue in the next subsection, much of the uncertainty arising from having a time dependency in $\opOpRfsdec$ leads to having a volatility that can be eliminated by just considering a time independent base site. This also eliminates the problem of having to work with Brownian motions. Thus instead of considering an element $\cF^f$ of $\opOpRfsdec$ we consider instead $\Omega^+$ made into a site by taking a Grothendieck topology very much like the one induced by the probability $P$ in the previous section. Then random variables will be upgraded to being sections of some sheaf over $\Omega^+$. Events being simplicial sets we call a morphism $\omega' \rarr \omega$ a structural morphism if it is a monomorphism of simplicial sets. We define coverings to be collections of such monomorphisms. Let $\tau$ be the collection of such coverings.

\begin{HONO_top}
$\tau$ defines a Grothendieck topology on $\Omega^+$ making it into a site.
\end{HONO_top}
\begin{proof}
If $\omega' \rarr \omega$ is an isomorphism then $\{ \omega' \rarr \omega \}$ is trivially a covering. If $\{ \omega_i \rarr \omega \}$ is a covering, $\gamma \rarr \omega$ any morphism, then as before $\omega_i \times_{\omega} \gamma$ is a simplicial set, so we have a morphism $\omega_i \times_{\omega} \gamma \rarr \gamma$ in $\Omega^+$ and such a morphism is a structural morphism as well, and this for all $i \in I$ so $\{ \omega_i \times_{\omega} \gamma \rarr \gamma \}$ is a covering of $\gamma$. Finally if $\{\omega_i \rarr \omega \}$ is a covering, $\{ \omega_{ij} \rarr \omega_i \}$ is a covering, then so is $\{\omega_{ij} \rarr \omega_i \rarr \omega \}$.
\end{proof}

\section{Construction of the sheaf of Brownian motions} \label{Constr}
The problem we are now facing is that of determining in what space, loosely speaking, are Brownian motions $X$ on $(\Omega, \cF, P)$ defined. As pointed out above it is necessary to consider $(\cC, \tau)$ instead with $\cC = \cF^f \in \opOpRfsdec$ and $\tau$ being given by the topology induced by $P$.\\

\subsection{Differentials}

In a first time observe that it is typical in Quantitative Finance to consider the differential of a stock price $S$ being given by $dS = \alpha S dt + \sigma S dW$ (\cite{S}), where $W$ is a Brownian motion. This means $S$ itself depends on $W$ in addition to depending on time. Note that this notation is actually misleading. The differential of a Brownian motion is ill-defined and is instead defined by considering for a given interval a subdivision $\Pi = \{t_i\}_{1 \leq i \leq n}$ of that interval. We then define as in \cite{S}:
\beq
dW = \lim_{\| \Pi \| \rarr 0} \sum_i \Delta_i W
\eeq
where $\Delta_i W = W(t_{i+1}) - W(t_i)$. Writing $\Delta$ without index and $W(t_{i+1})$ as $W_{i+1}$ observe that this implies:
\begin{align}
\Delta X \cdot \Delta Y &= (X_{i+1} - X_i)(Y_{i+1} - Y_i) \nonumber \\
X \Delta Y &= X_i(Y_{i+1} - Y_i) \nonumber \\
\Delta X \cdot Y &= (X_{i+1} - X_i)Y_i
\end{align}
summing those lines we end up with $X_{i+1}Y_{i+1} - X_i Y_i = \Delta(XY)$. Summing over the index $i$ we get $d(XY) = XdY + dX \cdot Y + dX \cdot dY$ as it is usual in stochastic calculus (\cite{S}). This means in particular this differential on Brownian motions is not a derivation in the usual algebraic sense, so henceforth we will denote it by $\delta$ instead of using $d$. Accordingly we will write:
\beq
DS = \alpha S dt + \sigma S \delta W
\eeq
An immediate question is whether we have $\delta ^2 = 0$. We will answer this question shortly.\\

Note also that when one develops Ito's formula (\cite{I}) for differentiating the function of a Brownian motion, one writes:
\beq
df(t,W) = \partial_1 f dt + \partial_2 f dW + \frac{1}{2} \partial_2 ^2 f (dW)^2
\eeq
so one assumes we can have a Taylor expansion of $f$ for which loosely speaking $dW dt = 0$ and $(dW)^2 = dt$. This subsumes that for the purposes of doing a Taylor expansion $\delta W$ is considered to be a finite difference, which it is by definition. \\

At this point we can perform a computation which will be needed later. For $X$ a random variable whose differential is given by $DX = \alpha X dt + \sigma X \delta W$ then using the above differential rule we have:
\begin{align}
D \log X &= \frac{1}{X} dX - \frac{1}{2} \frac{1}{X^2}(dX)^2 \nonumber \\
&=\alpha dt + \sigma \delta W - \frac{1}{2X^2} \sigma^2 X^2 dt \nonumber \\
&= (\alpha - \frac{1}{2}\sigma^2)dt + \sigma \delta W
\end{align}

\subsection{Construction of a sheaf of Brownian motions}

Recall that $\sigma$ stands for the volatility of $S$. If $W$ itself is a function of events, it has its own volatility. In probability events are well-defined. Not so in Finance. Either an event $\omega$ is trivial, in which case its complement $\omega^c$ is infinitely complex, or $\omega$ itself is infinitely complex, hence a function of events is bound to be highly volatile, a reflection that in practice events are not fully defined and there is consequently a large room for error which we can accomodate by seeing $W$ itself as having a volatility, which is passed on to the price of a stock, which itself adds some volatility of its own. Indeed as $dt \rarr 0$ one should be able to resolve events, hence resolve the uncertainty on $W$, but saying $W$ is a Brownian motion does not allow that. As a matter of fact the curve for $W$ cannot be resolved (\cite{L}). Thus in developing a theory of Brownian motions, this is one salient feature we have to incorporate in our formalism.\\

As just mentioned $W$ depends on events so $W$ is really a section of some space over $(\cC, \tau)$ which we will argue is a sheaf in a sense to be precised, and we regard $\delta \omega$ to be a deformation of $\omega$ along a path in the base. We would like to understand what $\delta W$ means. For $X$ and $Y$ two Brownian motions we have the well-known result (\cite{S}):
\beq
\delta(XY) = X \delta Y + \delta X . Y + \delta X \delta Y
\eeq
A first observation is that we have elements of the form $X \delta Y$ with $X$ and $\delta Y$ objects of a different nature so we should really write that as $X \otimes \delta Y$. That rule just given is typically taken as a formula when actually it is just a consequence of the definition of $\delta X$ for $X$ a Brownian motion. Recall that $\delta X$ is defined as the limit of $\Delta X$ when the subdivision size goes to zero, and $\Delta X$ we regard as the boundary $\partial X$ of $X$ where $\partial$ is the boundary operator on simplices. Thus $\delta X$ carries some homological information relative to the space of random variables. Further if we denote by $\cH$ the space in which Brownian motions are defined then we should write $XY = X \otimes Y$ with $\otimes$ a tensor product on $\cH$. Without loss of generality we can put a structure of a symmetrical monoidal category on $\cH$.\\

Collecting things together we first define a constant sheaf $\cW_t = \cR$ valued in $\mR$ over $\cF^f_t$ for all $t$, the collection of which
\beq 
 \cW = \{\cW_t \} \in \text{Map}(\mR, \text{Sh}(\opOpsdec)) \nonumber
\eeq 
is what we call a filtered sheaf. Now $\cW$ is not a sheaf as for $\omega$ fixed, $t < t'$, there is no unique object $W_{t'}(\omega)$ with a morphism $W_{t'}(\omega) \rarr W_t(\omega)$. What we have instead is a cone at each point $W_t(\omega) \in \cW_t(\omega)$, a cone with apex at $W_t(\omega)$, based at $t'$ in $\cW_{t'}(\omega)$ in which possible sections of $\cW_{t'}(\omega)$ originating from $\cW_t(\omega)$ are valued, the base of that cone being dependent on the increment $t' - t$ as well as the volatility of $W$. What that means is possible movements in $W_{t'}$ given $W_t(\omega)$ are found in the base of that cone. Thus in $\cW$ instead of having morphisms in $t$ we have transversal cones. This prompts us to generalize this formalism as follows.\\

\begin{ShBrown}
For $\cC \in \mapRcat$, $\tau$ a filtered Grothendieck topology on $\cC$, $\cD$ a fixed category, $\cW_t$ a sheaf in $\cD$ over $\cC(t)$ for all $t$, assembling into a filtered sheaf $\cW = \{ \cW_t \}$ over $(\cC, \tau)$, if we have transversal cones in time defined for all $W \in \Gamma(\cC, \cW)$ then such a filtered sheaf $\cW$ we call a sheaf of Brownian motions.
\end{ShBrown}

On the sections of $\cW$ we define a symmetric monoidal category structure. Sections of $\cW$ are regarded as geometric objects that carry some homological information. In particular we have a boundary map $\delta$ on $\Gamma(\cC, \cW)$. The set of points $W_t(\omega) \in \cW_t(\omega)$ inherit the monoidal structure from that of $\Gamma(\cC, \cW)$, but are only used in the definition of the homology of elements of such a space. Transversally, for $\omega$ fixed, $t$ fixed, $t'$ chosen, $\delta W$ at $(t, \omega)$ with respect to $t'$ is $W_{t'}(\omega) - W_t(\omega)$. We can also defined a boundary map in the shaves $\cW_t$ as follows: for $\omega$ fixed, $\omega \rarr \omega'$ a morphism, $\delta \omega$ along that morphism is defined to be $\omega' - \omega$, which induces $\delta W_t = W_t(\omega') - W_t(\omega)$. Note that this provides a generalization of $\delta X$ that is typically not considered in Finance. We can make the notion of transversal homology more precise as follows: recall that for a given category $\cC$ we have the Yoneda embedding $j: \cC \rarr \text{Set}^{\cC^{op}}, C \mapsto \text{Hom}_{\cC}( \cdot, C)$. For us $\cC$ is played by the $\cC$-diagram in a category $\cD$ formed by a functor $F: \cC \rarr \cD$, which we regard as a simplicial set, in which case denoting this simplicial set by $F(\cC)$ we still have an $\infty$-generalization of the Yoneda embedding: $F(\cC) \rarr \cS^{F(\cC)^{op}}$, $\cS$ the $\infty$-category of spaces (\cite{Lu2}). Then for $c, c' \in \cC$, $F(c') - F(c) \mapsto \text{Hom}_{F(\cC)}( \cdot, F(c')) - \text{Hom}_{F(\cC)}( \cdot, F(c))$ under the Yoneda embedding. Those presheaves are objects in $\text{Fun}(F(\cC)^{op}, \cS)$, an $\infty$-category since $\cS$ is an $\infty$-category, hence we can regard $F(c')-F(c)$ as the boundary of a morphism $F(c) \rarr F(c')$ in $\text{Fun}(F(\cC)^{op}, \cS)$, which we take to be $F$, which leads us to the following definition:
\begin{Deform}
For $\cC$ and $\cD$ two categories, $F: \cC \rarr \cD$ a functor, we define the $q$-boundary $\delta F$ of $F$ at $c \in \text{Ob}(\cC)$ along a morphism $c \rarr c'$ in $\cC$ to be defined by $\delta F(c) = F(c') - F(c)$.
\end{Deform}

Note that if this definition gives a meaning to having a difference $F(c')-F(c)$, it does not explain how such a difference should be used in practice. For a difference such as $W_{t'}(\omega) - W_t(\omega)$ Brownian motions being real valued the difference is a real number. For more general categories $\cD$ it is really contextually that one can have a working definition of $F(c')- F(c)$.\\

For a fixed event $\omega$, $W_t(\omega)$ as a function of the index $t$ will display the characteristics of a Brownian motion. Thus letting $W_t(\omega) = W(t, \omega)$ for a Brownian motion $W$, we picture a cone at $W(t, \omega)$ in the $(t, W)$-plane opening to the right whose opening is dependent upon $W$'s volatility. Possible elements $W_{t'}(\omega)$ for $t' > t$ are in that cone. Observe that for $t$ fixed there are morphisms between events in $\cF^f_t$, which necessarily means some time has elapsed when those sets are assembled, while $t$ itself is fixed. This is no contradiction. Recall that $t$ is an index, that we have a bundle $\mR^f \rarr \mR $, and it is in the fiber of this map that we have morphisms between events. Thus at $t$ we have an $(0,1]$-fiber in time in which morphisms in $\cF^f_t$ occur. This also means that for two indexes $t$ and $t'$, $t < t'$, $\omega$ fixed, the transition $W(t, \omega) \rarr W(t', \omega)$ is uncertain as we are missing all the morphisms between events in the fibers over points $s$, $t < s < t'$, whence the uncertainty in defining $W$ and $\cW_t$, which gives the Brownian motion picture.\\

\subsection{Tropical Geometry realization}
As a digression we consider the tropical realization (\cite{M}) of stochastic differential equations. If we have $DX = \alpha X dt + \sigma X \delta W$, that corresponds to writing $D \log X = (\alpha - \frac{1}{2}\sigma^2 )dt + \sigma \delta W$. Having a log differential invites us to consider the tropicalization of this equation:
\beq
D \log X = ``(\alpha - \frac{1}{2}\sigma^2)dt + \sigma \delta W" = \max \{ \alpha - \frac{1}{2}\sigma^2 + dt , \sigma + \delta W \}
\eeq
\begin{cdga}
For $\cA$ a $\mathbb{C}\text{-dga}$, $\alpha = \sum \alpha_i \nu^i \in \cA$, $\alpha_i \in \mathbb{C}$, we define an augmentation homomorphism:
\begin{align}
\varepsilon: \cA &\rarr \mathbb{C} \nonumber \\
\alpha &\dashrightarrow \sum \alpha_i
\end{align}
Then if $\beta = \sum \beta_i \nu^i \in \cA$ we define:
\beq
\text{max}\{\alpha, \beta \} = \text{max} \{ \varepsilon(\alpha), \varepsilon(\beta) \}
\eeq
\end{cdga}
With this definition we have:
\beq
D \log X = \text{max} \{ \alpha - \frac{1}{2} \sigma^2 + 1 , \sigma +1 \} = \text{max} \{ \alpha - \frac{1}{2} \sigma^2 , \sigma \}
\eeq

If we work on $\cF^f \in \opOpRfsdec$, then this is the problem we are facing. As an interesting note observe that $\alpha$ and $\sigma$ are put on the same footing in the tropical picture, which is not anything new as we have the market price of risk variable $\theta = (\alpha - r)/\sigma$ (\cite{S}). If we work with $\Omega^+$ however, there is no time dependency, thus no volatility, no Brownian motion and we just have a linearization $D \log X = \Delta \log X$ for deformations of events, something we will develop in the last section.\\

In a categorical setting we need to define the logarithm of objects and morphisms. For $X$ a random variable on $(\cF^f, \tau_P) \in \opOpRfsdec$ we formally define:
\beq
\exp(X) = \sum_{n \geq 0} \frac{1}{n!}\otimes^n X
\eeq
where we would have a symmetrical monoidal structure on a filtered sheaf of random variables $\cG$ on $(\cF^f, \tau_P)$, so we would have a functor
\begin{align}
\exp: \cG &\rarr \overline{\cG}  \nonumber \\
X &\mapsto \exp(X)
\end{align}
where $\overline{\cG}$ is the graded completion of the symmetric monoidal category associated with $\cG$. Then the formal inverse to such a functor would be a log morphism:
\begin{align}
\cG &\leftarrow \overline{\cG} : \log \nonumber \\
X &\dashleftarrow \exp(X)
\end{align}
Formally one has $\log X = \sum_{n>0}(-1)^n/n \otimes^n X$. On morphisms if $X$ and $Y$ are objects of $\cG$ and $f$ is a morphism from $X$ to $Y$ then $\exp(f) = \sum 1/n! \otimes^n f$ and likewise for $\log(f)$.

\section{Of the use of $\Omega^+$ to resolve Wiener processes}

Note that if we have an event $\omega$ in a $\sigma$-algebra over $\Omega$ then $\omega^c = \Omega - \omega$ is also in that $\sigma$-algebra. If $\omega$ is simple enough, to have a function defined at $\omega^c$ one would have to resort to an approximation. Further having filtered objects as we have seen introduces a further uncertainty in defining transition maps, so instead of working with filtered $\sigma$-algebras and probabilities, one should work with categories and structural morphisms as done with $\Omega^+$ in \ref{TIF}.\\

One can retain a notion of time in this time independent formalism by introducing roof diagrams such as the one below where $\lhd(A)$, which we call the forward cone at $A$, is defined to be the collection of events $B$ that can be reached from $A$ via a morphism $A \rarr B$ in $\Omega^+$. Thus we consider a category $\cE$ with $\text{Ob}(\cE) = \text{Ob}(\Omega^+)$ and morphisms in $\cE$ are roof diagrams such as the one below with base an element of $\text{Mor}(\Omega^+)$.

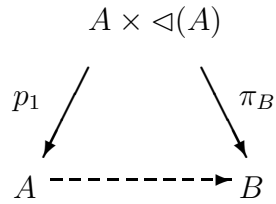
\beq
\setlength{\unitlength}{0.5cm}
\begin{picture}(6,5)(0,0)
\thicklines
\put(-1,-0.5){$A$}
\put(5,-0.5){$B$}
\put(1,4){$A \times \lhd(A)$}
\multiput(0,0)(0.5,0){10}{\line(1,0){0.3}}
\put(4.5,0){\vector(1,0){0.3}}
\put(1,3){\vector(-1,-2){1.2}}
\put(4,3){\vector(1,-2){1.2}}
\put(-1,2){$p_1$}
\put(5,2){$\pi_B$}
\end{picture}
\nonumber
\eeq

Identity:

\beq
\setlength{\unitlength}{0.5cm}
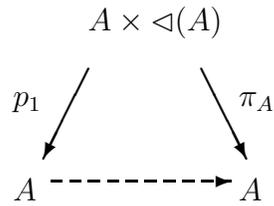
\begin{picture}(6,5)(0,0)
\thicklines
\put(-1,-0.5){$A$}
\put(5,-0.5){$A$}
\put(1,4){$A \times \lhd(A)$}
\multiput(0,0)(0.5,0){10}{\line(1,0){0.3}}
\put(4.5,0){\vector(1,0){0.3}}
\put(1,3){\vector(-1,-2){1.2}}
\put(4,3){\vector(1,-2){1.2}}
\put(-1,2){$p_1$}
\put(5,2){$\pi_A$}
\end{picture}
\nonumber
\eeq

Composition:

\beq
\setlength{\unitlength}{0.5cm}
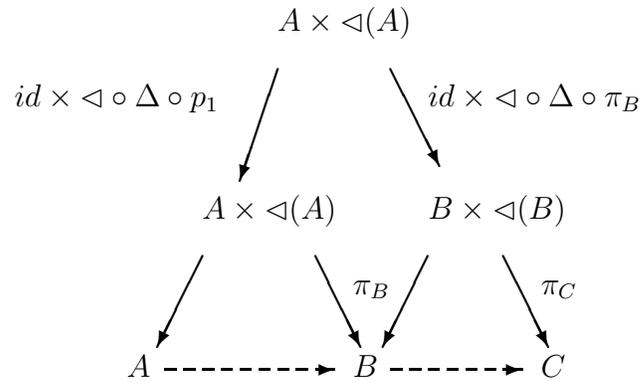
\begin{picture}(11,10)(0,0)
\thicklines
\put(-1,-0.2){$A$}
\put(5,-0.2){$B$}
\put(10,-0.2){$C$}
\put(1,4){$A \times \lhd(A)$}
\put(7,4){$B \times \lhd(B)$}
\multiput(0,0)(0.5,0){8}{\line(1,0){0.3}}
\put(4,0){\vector(1,0){0.5}}
\multiput(6,0)(0.5,0){7}{\line(1,0){0.3}}
\put(9,0){\vector(1,0){0.5}}
\put(1,3){\vector(-1,-2){1.2}}
\put(4,3){\vector(1,-2){1.2}}
\put(7,3){\vector(-1,-2){1.2}}
\put(9,3){\vector(1,-2){1.2}}
\put(5,2){$\pi_B$}
\put(10,2){$\pi_C$}
\put(3,9){$A \times \lhd(A)$}
\put(3,8){\vector(-1,-3){1}}
\put(-4,7){$id \times \lhd \circ \Delta \circ p_1$}
\put(6,8){\vector(1,-2){1.3}}
\put(7,7){$id \times \lhd \circ \Delta \circ \pi_B$}
\end{picture}
\nonumber
\eeq

where $\Delta$ is the diagonal map. With this composition being defined we can verify that the identity defined above functions as an identity:

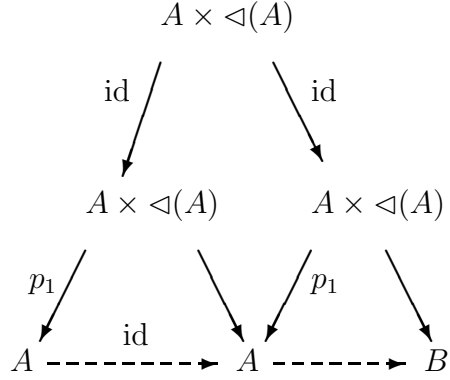
\beq
\setlength{\unitlength}{0.5cm}
\begin{picture}(11,10)(0,0)
\thicklines
\put(-1,-0.2){$A$}
\put(5,-0.2){$A$}
\put(10,-0.2){$B$}
\put(1,4){$A \times \lhd(A)$}
\put(-0.5,2){$p_1$}
\put(7,4){$A \times \lhd(A)$}
\put(7,2){$p_1$}
\multiput(0,0)(0.5,0){8}{\line(1,0){0.3}}
\put(2,0.5){$\text{id}$}
\put(4,0){\vector(1,0){0.5}}
\multiput(6,0)(0.5,0){7}{\line(1,0){0.3}}
\put(9,0){\vector(1,0){0.5}}
\put(1,3){\vector(-1,-2){1.2}}
\put(4,3){\vector(1,-2){1.2}}
\put(7,3){\vector(-1,-2){1.2}}
\put(9,3){\vector(1,-2){1.2}}
\put(3,9){$A \times \lhd(A)$}
\put(1.5,7){$\text{id}$}
\put(3,8){\vector(-1,-3){1}}
\put(6,8){\vector(1,-2){1.3}}
\put(7,7){$\text{id}$}
\end{picture}
\nonumber
\eeq

which simplifies as:

\beq
\setlength{\unitlength}{0.5cm}
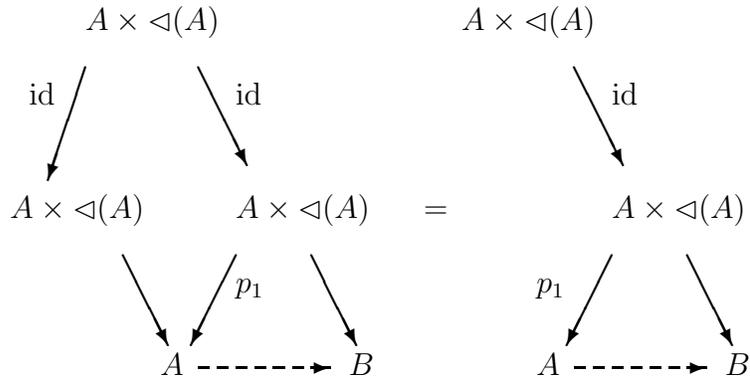
\begin{picture}(20,10)(0,0)
\thicklines
\put(5,-0.2){$A$}
\put(10,-0.2){$B$}
\put(1,4){$A \times \lhd(A)$}
\put(7,4){$A \times \lhd(A)$}
\put(7,2){$p_1$}
\multiput(6,0)(0.5,0){7}{\line(1,0){0.3}}
\put(9,0){\vector(1,0){0.5}}
\put(4,3){\vector(1,-2){1.2}}
\put(7,3){\vector(-1,-2){1.2}}
\put(9,3){\vector(1,-2){1.2}}
\put(3,9){$A \times \lhd(A)$}
\put(1.5,7){$\text{id}$}
\put(3,8){\vector(-1,-3){1}}
\put(6,8){\vector(1,-2){1.3}}
\put(7,7){$\text{id}$}
\put(12,4){$=$}
\put(15,-0.2){$A$}
\put(20,-0.2){$B$}
\put(17,4){$A \times \lhd(A)$}
\put(15,2){$p_1$}
\put(17,3){\vector(-1,-2){1.2}}
\multiput(16,0)(0.5,0){7}{\line(1,0){0.3}}
\put(19,0){\vector(1,0){0.5}}
\put(19,3){\vector(1,-2){1.2}}
\put(13,9){$A \times \lhd(A)$}
\put(16,8){\vector(1,-2){1.3}}
\put(17,7){$\text{id}$}
\end{picture}
\nonumber
\eeq

which further simplifies as:

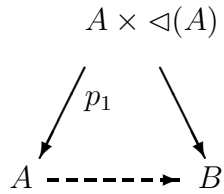
\beq
\setlength{\unitlength}{0.5cm}
\begin{picture}(5,5)(0,0)
\thicklines
\put(0,-0.2){$A$}
\put(5,-0.2){$B$}
\put(2,4){$A \times \lhd(A)$}
\put(2,2){$p_1$}
\put(2,3){\vector(-1,-2){1.2}}
\multiput(1,0)(0.5,0){7}{\line(1,0){0.3}}
\put(4,0){\vector(1,0){0.5}}
\put(4,3){\vector(1,-2){1.2}}
\end{picture}
\nonumber
\eeq

thereby showing that the identity as we defined it acts as a left identity and likewise one would show that it also acts as a right identity.\\

Associativity:

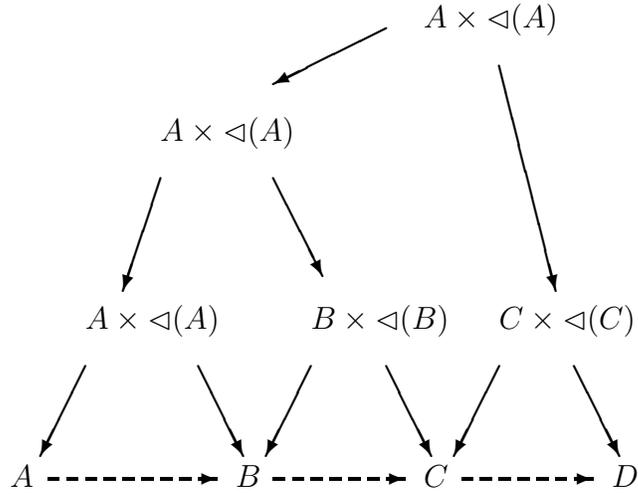
\beq
\setlength{\unitlength}{0.5cm}
\begin{picture}(16,14)(0,0)
\thicklines
\put(-1,-0.2){$A$}
\put(5,-0.2){$B$}
\put(10,-0.2){$C$}
\put(15,-0.2){$D$}
\put(1,4){$A \times \lhd(A)$}
\put(7,4){$B \times \lhd(B)$}
\put(12,4){$C \times \lhd(C)$}
\multiput(0,0)(0.5,0){8}{\line(1,0){0.3}}
\put(4,0){\vector(1,0){0.5}}
\multiput(6,0)(0.5,0){7}{\line(1,0){0.3}}
\put(9,0){\vector(1,0){0.5}}
\multiput(11,0)(0.5,0){7}{\line(1,0){0.3}}
\put(14,0){\vector(1,0){0.5}}
\put(1,3){\vector(-1,-2){1.2}}
\put(4,3){\vector(1,-2){1.2}}
\put(7,3){\vector(-1,-2){1.2}}
\put(9,3){\vector(1,-2){1.2}}
\put(12,3){\vector(-1,-2){1.2}}
\put(14,3){\vector(1,-2){1.2}}
\put(3,9){$A \times \lhd(A)$}
\put(3,8){\vector(-1,-3){1}}
\put(6,8){\vector(1,-2){1.3}}
\put(10,12){$A \times \lhd(A)$}
\put(9,12){\vector(-2,-1){3}}
\put(12,11){\vector(1,-4){1.5}}
\end{picture}
\nonumber
\eeq

and for the associated composition:

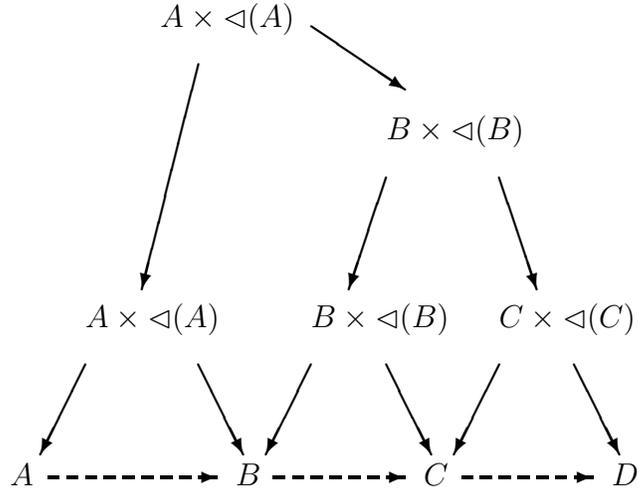
\beq
\setlength{\unitlength}{0.5cm}
\begin{picture}(16,14)(0,0)
\thicklines
\put(-1,-0.2){$A$}
\put(5,-0.2){$B$}
\put(10,-0.2){$C$}
\put(15,-0.2){$D$}
\put(1,4){$A \times \lhd(A)$}
\put(7,4){$B \times \lhd(B)$}
\put(12,4){$C \times \lhd(C)$}
\multiput(0,0)(0.5,0){8}{\line(1,0){0.3}}
\put(4,0){\vector(1,0){0.5}}
\multiput(6,0)(0.5,0){7}{\line(1,0){0.3}}
\put(9,0){\vector(1,0){0.5}}
\multiput(11,0)(0.5,0){7}{\line(1,0){0.3}}
\put(14,0){\vector(1,0){0.5}}
\put(1,3){\vector(-1,-2){1.2}}
\put(4,3){\vector(1,-2){1.2}}
\put(7,3){\vector(-1,-2){1.2}}
\put(9,3){\vector(1,-2){1.2}}
\put(12,3){\vector(-1,-2){1.2}}
\put(14,3){\vector(1,-2){1.2}}
\put(9,9){$B \times \lhd(B)$}
\put(9,8){\vector(-1,-3){1}}
\put(12,8){\vector(1,-3){1}}
\put(3,12){$A \times \lhd(A)$}
\put(4,11){\vector(-1,-4){1.5}}
\put(7,12){\vector(3,-2){2.5}}
\end{picture}
\nonumber
\eeq

We have associativity by virtue of the equality $\pi_{B \rarr C} \circ \pi_{A \rarr B} = \pi_{A \rarr C}$. We put on $\cE$ a Grothendieck topology generated by roof diagrams with base structural morphisms - which we call structural roof diagrams - as introduced in section \ref{TIF}. We verify that this does indeed give such a Grothendieck topology. For coverings we take collections of structural roof diagrams $\{\omega' \rarr \omega \}$. If we have an isomorphism $\omega' \rarr \omega$ then $\{\omega' \rarr \omega \}$ is a covering. If $\{\omega_i \rarr \omega \}$ is a covering and $\gamma \rarr \omega$ is any morphism for all $i$ the projections $\omega_i \times_{\omega} \gamma \rarr \gamma$ are structural morphisms, so all such morphisms form the base of a covering of $\gamma$. Finally if $\{\omega_i \rarr \omega \}$ is a covering, and so are the $\{ \omega_{ij} \rarr  \omega_i \}$ then so is the composition $\omega_{ij} \rarr \omega_i \rarr \omega$, so the collection of all such compositions is a covering of $\omega$. Denoting by $\tau$ the Grothendieck topology thus defined, $(\cE, \tau)$ is a site.\\

Finally we consider the category of sheaves $\cW: (\cE, \tau)^{op} \rarr \text{Cat}^{\otimes}_{\infty}$ of symmetric monoidal $\infty$-categories. Since there is no time dependence, no probability, no $\sigma$-algebra, all sources of uncertainties are eliminated.\\

Finally for $X$ a section of such a sheaf $\cW$, $\omega \in \cE$, we let $\omega^o$ be the set of morphisms $\omega \rarr \omega'$ in a same connected component $C$ $\omega$ is in for which there is no $\omega'' \in C$ such that $\omega \rarr \omega'' \rarr \omega$. Then for $\psi: \omega \rarr \omega' \in \omega^o$, we define:
\beq
D_{\psi}X(\omega) = X(\omega') - X(\omega)
\eeq
If one can define quotients in $\cW$ then one would define this difference as being $X(\omega') / X(\omega)$. Note that if that solves the problem of having to solve a log differential equation in Stochastic calculus, this is not very illuminating however as it is quite a simple equation. What is most important perhaps is how $X$ behaves in between $\omega$ and $\omega'$ relative to other sections of $\cW$ between the same objects, and this from a homological perspective, so one could impose that our sheaves be sheaves of stable symmetric monoidal $\infty$-categories (\cite{Lu3}) and study their sheaf cohomology. To be complete one should also study such a sheaf as an object of the topos of symmetric monoidal $\infty$-categories (\cite{Lu2}) on $\cE$.

\end{document}